\documentclass{amsart}
\usepackage{amsfonts,amssymb,amscd,amsmath,enumerate,verbatim}
\usepackage[all]{xy}
\SelectTips{eu}{}
\xyoption{curve}

\newcommand\BN{{\mathbb N}}
\newcommand\BZ{{\mathbb Z}}

\newcommand\ck{{\mathsf K}}
\newcommand\cs{{\mathsf S}}
\newcommand\ct{{\mathsf T}}
\newcommand\cd{{\mathsf D}}

\newcommand\ssc{{\mathsf c}}

\newcommand{\xra}{\xrightarrow}
\newcommand{\lra}{\longrightarrow}

\newcommand{\les}{{\scriptscriptstyle\leqslant}}

\newcommand{\sls}{{\scriptscriptstyle<}}

\newcommand{\shift}{{\sf\Sigma}}
\newcommand\col{\colon}

\newcommand\dd{\partial}
\newcommand{\hh}[1]{\operatorname{H}(#1)}

\newcommand{\HH}[2]{\operatorname{H}_{#1}(#2)}

\newcommand{\Hom}[3]{\operatorname{Hom}_{#1}(#2,#3)}

\renewcommand\ker{\operatorname{Ker}}

\newcommand\idmap{\operatorname{id}}
\newcommand\fd{\operatorname{fd}}
\newcommand\hfd{\operatorname{gr-fd}}
\newcommand\id{\operatorname{id}}
\newcommand\hid{\operatorname{gr-id}}

\newcommand\pd{\operatorname{pd}}
\newcommand\hpd{\operatorname{gr-pd}}

\newcommand{\dcat}[1][R]{{\mathsf{D}(#1)}}
\newcommand{\dcatc}[1][R]{{\mathsf{D}^{\mathsf c}(#1)}}
\newcommand{\dcatf}[1][R]{{\mathsf{D}^{\mathsf f}_{\mathsf b}(#1)}}
\newcommand{\dcatfp}[1][R]{{\mathsf{D}^{\mathsf{fp}}_{\mathsf b}(#1)}}

\newcommand{\kac}[1][R]{{\mathsf{K}_{\mathsf{ac}}(\operatorname{Inj}\,#1)}}
\newcommand{\kinj}[1][R]{{\mathsf{K}(\operatorname{Inj}\,#1)}}

\newcommand{\kinjc}[1][R]{{\mathsf{K}^{\mathsf c}(\operatorname{Inj}\,#1)}}
\newcommand{\kpac}[1][R]{{\mathsf{K}_{\mathsf{ac}}(\operatorname{Proj}\,#1)}}
\newcommand{\kproj}[1][R]{{\mathsf{K}(\operatorname{Proj}\,#1)}}
\newcommand{\kprojc}[1][R]{{\mathsf{K}^{\mathsf c}(\operatorname{Proj}\,#1)}}

\newcommand{\op}[1]{{#1}^{\mathrm{op}}}

\newcommand{\eps}{{\varepsilon}}

\theoremstyle{plain}
\newtheorem{theorem}{Theorem}[section]

{\Alph{theorem}}
\newtheorem{proposition}[theorem]{Proposition}
\newtheorem{lemma}[theorem]{Lemma}

\newtheorem{corollary}[theorem]{Corollary}

\theoremstyle{definition}
\newtheorem*{ack}{Acknowledgments}

\newtheorem{chunk}[theorem]{}

\theoremstyle{remark}

\newtheorem{remark}[theorem]{Remark}
\numberwithin{equation}{theorem}

\begin{document}
\title[Homological dimensions]
{Homological dimensions and regular rings}


\author[A.~Iacob]{Alina Iacob}
\address{Department of Mathematical Sciences,
 Georgia Southern University,
     Statesboro, GA 30460, USA}
\email{aiacob@georgiasouthern.edu}

\author[S.~B.~Iyengar]{Srikanth B.~Iyengar}
\address{Department of Mathematics, University of Nebraska, Lincoln, NE 68588, U.S.A.}
\email{iyengar@math.unl.edu}

\thanks {S.B.I. was partly supported by NSF grants DMS 0602498 and DMS 0903493}

  \begin{abstract}
A question of Avramov and Foxby concerning injective dimension of complexes is settled in the affirmative for the class of noetherian rings. A key step in the proof is to recast the problem on hand into one about the homotopy category of complexes of injective modules. Analogous results for flat dimension and projective dimension are also established.
  \end{abstract}

\keywords{homological dimension, regular ring, semi-injective, semi-projective}
\subjclass[2000]{13D25, 13D02}

\maketitle

\section{Introduction}
There is a well-established notion of injective dimension for
modules over a ring, based on a natural construction of injective
resolutions for modules. However, in extending it from modules to
complexes one can choose from various extensions of the notion of a
resolution, and these yield potentially different notions of
injective dimension. The resulting concepts were defined by Avramov
and Foxby in \cite{Avramov/Foxby:1991a} where they proved that they
yield the same invariants when the ring has finite global dimension,
and asked if the converse statement holds. In this article we settle
this question for certain classes of rings. To illustrate the issues
involved, in the remainder of the Introduction we focus on the case
of modules.

Let $R$ be a ring and $M$ an $R$-module; in what follows `module'
means `left module' and properties considered are with respect to
the left structures, unless stated otherwise. Recall that a
classical injective resolution of $M$ is a complex $I$ of injective
$R$-modules with $I_j=0$ for $j > 0$ and $\hh I\cong M$, and that
the injective dimension of $M$, denoted $\id_{R}M$, is the infimum
of those integers $n$ such that $M$ admits an injective resolution
$I$ such that $I_j=0$ for $j<-n$. We consider a variant of this
notion where $I$ is no longer limited to non-positive degrees:
\[
\inf
\left\{n\in\BZ\,  \left|
\begin{gathered}
\text{$M\cong \hh I$ with $I$ a complex of}\\
\text{injectives with $I_{j}=0$ for $j<-n$.}
\end{gathered}
\right. \right\}
\]
Following \cite{Avramov/Foxby/Halperin:2008a} we denote\footnote{In~\cite{Avramov/Foxby:1991a} this number is denoted $\#\!-\!\mathrm{id}_R M$.} this number $\hid_{R}M$. Evidently, there is an inequality
\[
\hid_{R}M \leq \id_{R}M\,.
\]
Avramov and Foxby~\cite[3.5]{Avramov/Foxby:1991a} proved that equality holds whenever the ring $R$ has finite global dimension, and asked---see  Question 3.8 in \emph{op.~cit.}---if the converse is true.
We answer this question in the affirmative in the class of noetherian rings:

\begin{theorem}
\label{ithm:inj}
When $R$ is noetherian the following conditions are equivalent.
\begin{enumerate}[\quad\rm(1)]
\item $R$ is regular.
\item $\hid_{R}M=\id_{R}M $ for every module $M$.
\end{enumerate}
\end{theorem}

Recall that the ring $R$ is said to be \emph{regular} if every ideal has a finite resolution by finitely generated projective modules; equivalently, if $R$ is noetherian and the projective dimension of each finitely generated module is finite. There exist regular rings of infinite global dimension (see~\cite[Appendix, Example 1]{Nagata:1962a}) so (1)$\implies$(2) above strengthens, for noetherian rings, the result of Avramov and Foxby~\cite[3.5]{Avramov/Foxby:1991a}.

Theorem~\ref{ithm:inj} follows from Proposition~\ref{id=hid} and Theorem~\ref{thm:inj}. A key step in the proof is to recast the statement about modules in terms of properties of the homotopy category of complexes of injective modules. Results of Krause~\cite{Krause:2005a} on this homotopy category are then invoked to complete the argument.

The next result is an analogue of Theorem~\ref{ithm:inj} for projective and flat modules; it answers another part of \cite[3.8]{Avramov/Foxby:1991a}. In its statement $\pd_{R}M$ and $\fd_{R}M$ are the projective dimension and flat dimension of $M$ respectively, while $\hpd_{R}M$ and $\hfd_{R}M$ are the corresponding analogues of $\hid_{R}M$.  As usual, $R$ is said to be \emph{coherent} if finitely generated ideals in $R$ are finitely presented; equivalently, if finitely generated submodules of free modules are finitely presented. We write `$\op R$' for the opposite ring of $R$, so `$\op R$ coherent' means that $R$ is right coherent.

\begin{theorem}
\label{ithm:prj}
When $\op R$ is coherent the following conditions are equivalent.
\begin{enumerate}[\quad\rm(1)]
\item Each bounded complex of finitely presented $\op R$-modules is perfect.
\item $\hpd_{R}M=\pd_{R}M $ for every module $M$.
\item $\hfd_{R}M=\fd_{R}M $ for every module $M$.
\end{enumerate}
\end{theorem}

Recall that a \emph{perfect complex} is one that is quasi-isomorphic
to a bounded complex of finitely generated projective modules. Note
that when $\op R$ is noetherian, (1) is equivalent to the condition
that the ring $\op R$ is regular.

The theorem above is contained in Proposition~\ref{pd=hpd} and Theorem~\ref{thm:prj}.  As for Theorem~\ref{ithm:inj}, the crucial idea  is to recast the conditions in terms of a homotopy category, but this time  the homotopy category of complexes of projective modules. We then apply results of J{\o}rgensen~\cite{Jorgensen:2005a} and Neeman~\cite{Neeman:2008a} to complete the proof.

\section{Complexes of injectives}
\label{Injectives}

Let $R$ be a ring, and let $\dcat$ denote the derived category of (left) $R$-modules; see Verdier~\cite{Verdier:1996a} for a construction of the derived category. We write $M\simeq N$ to indicate that $M$ and $N$ are \emph{quasi-isomorphic} complexes of $R$-modules, that is to say, they are isomorphic in $\dcat$. A morphism of complexes $M\to  N$ is a quasi-isomorphism if and only if its mapping cone, say $C$, is \emph{acyclic}, that is to say, $\hh C=0$ holds.

We recall some notions and results from \cite{Avramov/Foxby:1991a}, with terminology borrowed from \cite{Avramov/Foxby/Halperin:2008a}.

Let $I$ be a complex of $R$-modules. We say that $I$ is
\emph{graded-injective} if each $R$-module $I_{n}$ is injective;
equivalently, if the graded $R$-module underlying $I$ is injective
in the category of graded $R$-modules. The complex $I$ is
\emph{semi-injective} if it is graded-injective and whenever $\phi\col M\to N$ is a
quasi-isomorphism of complexes,  so is  $\Hom R{\phi}I\col \Hom
RNI\to\Hom RMI$. For instance, when $I$ is graded-injective and
$I_{n}=0$ for $n\gg 0$, it is semi-injective. For each complex $M$
there exists a quasi-isomorphism $M\to I$ with $I$ semi-injective;
such a morphism can also be chosen to be one-to-one. See
\cite{Avramov/Foxby/Halperin:2008a} for proofs of the assertions
above.

Give a class $\mathcal I$ of complexes of $R$-modules, consider the number
\[
\inf
\left\{n\in\BZ \,|\,
\text{$M\simeq I$ with $I\in\mathcal I$ and $I_{j} =0$ for $j<-n$.}
\right\}
\]
Taking for $\mathcal I$ the class of semi-injective complexes one gets the \emph{injective dimension} of $M$, denoted $\id_{R}M$. Taking for $\mathcal I$ the class of graded-injective complexes, and not only the semi-injective ones, gives rise to an invariant of $M$ that we denote $\hid_{R}M$. When $M$ is a module,
viewed as a complex with $M$ in degree zero and zero otherwise, these definitions yield the same invariants as those in the Introduction. This is because $M\simeq I$ holds if and only if $M\cong \hh I$ holds.

It is obvious from definitions that an inequality $\hid_{R}M\leq \id_{R}M$ holds for each complex $M$.
The result below describes some conditions under which equality holds for all $M$. Contractible complexes of injectives, appearing in condition (4), are the \emph{categorically injective} complexes of \cite{Avramov/Foxby/Halperin:2008a}.

\begin{proposition}
\label{id=hid}
Let $R$ be a ring. The following conditions are equivalent.
\begin{enumerate}[\quad\rm(1)]
\item An equality $\hid_{R}M=\id_{R}M $ holds for each complex $M$ of $R$-modules.
\item An equality $\hid_{R}M=\id_{R}M $ holds for each $R$-module $M$.
\item Each complex of injective $R$-modules is semi-injective.
\item Each acyclic complex of injective $R$-modules is contractible.
\end{enumerate}
\end{proposition}

\begin{proof}
Clearly, (3)$\implies$(1) and (1)$\implies$(2) hold.

(2)$\implies$(4) Let $I$ be an acyclic complex of injective $R$-modules. For each integer $i$ the inclusion $\ker(\dd_{i})\to \shift^{-i}I_{\les i}$ is a quasi-isomorphism. Thus, $\hid_{R} \ker(\dd_{i})=0$ and so the hypothesis entails $\ker(\dd_{i})$ is injective. Hence $I$ is contractible.

(4)$\implies$(3) Let $I$ be a complex of injectives and let $\iota\col I \to J$ be a semi-injective resolution of $I$. Since $\iota$ is a quasi-isomorphism, its mapping cone, say $C$, is acyclic; since $C$ is also a complex of injectives, the hypothesis yields that $C$ is contractible. Thus $\iota$ is a homotopy equivalence, and hence $I$ is itself semi-injective.
\end{proof}

Next we translate the equivalent conditions in the preceding
proposition to a condition concerning the homotopy category of
complexes of injective $R$-modules, which we denote $\kinj$. Its
objects are complexes of injective $R$-modules and its morphisms are
homotopy classes of morphisms of complexes; see \cite{Verdier:1996a}
for details and for a description of the triangulated structure
carried by $\kinj$. Let $\kac$ be the subcategory of $\kinj$
consisting of acyclic complexes; it is a triangulated subcategory.
There then exists a canonical localization functor
\[
Q\col \kinj \to \dcat\,.
\]
Its kernel is precisely $\kac$, so the next result is obvious.

\begin{lemma}
\label{ikac=0}
The functor $Q\col \kinj \to \dcat$ is an equivalence if and only if each acyclic complex of injective $R$-modules is contractible. \qed
\end{lemma}

Let $\ct$ be a triangulated category admitting all coproducts. An object $X$ in $\ct$ is \emph{compact} if $\Hom{\ct}X-$ commutes with all coproducts in $\ct$. In what follows $\ct^{\ssc}$ denotes the full subcategory of all compact objects in $\ct$. We say that $\ct$ is \emph{compactly generated} if the isomorphism classes of compact objects form a set, and the smallest triangulated subcategory containing $\ct^{\ssc}$ and closed under all coproducts is $\ct$ itself. The reader may refer to Neeman's book~\cite{Neeman:2001a} for a  discussion of these concepts.

A proof of the result below can be found, for example, in \cite[\S5.3]{Keller:1994a}, or \cite[2.2]{Neeman:1992b}.

\begin{chunk}
\label{dcat:cg}
The triangulated category $\dcat$ is compactly generated, and the objects of $\dcatc$ are the perfect complexes of $R$-modules.
\end{chunk}

We recall a well known, and not difficult to prove, test for equivalence of compactly generated categories; see, for example, \cite[4.5]{Benson/Iyengar/Krause:2008b}, or \cite[\S4.2]{Keller:1994a}.

\begin{chunk}
\label{cg:functor}
Let $F\col \cs\to \ct$ be an exact functor of compactly generated triangulated categories which is compatible with  coproducts. The functor $F$ is then an equivalence of categories if and only if it restricts to an equivalence of categories $\cs^{\ssc}\to \ct^{\ssc}$.
\end{chunk}

Next we focus on the case the ring $R$ is noetherian. In our context, this property is relevant because of the following result, due to Bass~\cite[1.1]{Bass:1962a}.

\begin{chunk}
\label{bass}
The ring $R$ is noetherian if and only if any (equivalently, any countable) direct sum of injective $R$-modules is injective.
\end{chunk}

When $R$ is noetherian we write $\dcatf$ for the full subcategory of $\dcat$ consisting of complexes $M$ such that $\HH iM$ is finitely generated for each $i$ and equal to zero when $|i|\gg0$ holds. Given \ref{dcat:cg}, the next remark is obvious.

\begin{chunk}
\label{dcat:regular}
A noetherian $R$ is regular if and only if $\dcatc=\dcatf$ holds.
\end{chunk}

The next result is due to Krause~\cite[2.3]{Krause:2005a}.

\begin{chunk}
\label{kinj:cg}
When $R$ is noetherian, $\kinj$ is compactly generated and the localization functor $Q$ induces an equivalence of categories
\[
Q\col \kinjc \xra{\ \sim \ } \dcatf\,.
\]
\end{chunk}

Proposition~\ref{id=hid} and the next result contain Theorem~\ref{ithm:inj} from the Introduction.

\begin{theorem}
\label{thm:inj}
Let $R$ be a ring. The following conditions are equivalent.
\begin{enumerate}[\quad\rm(1)]
\item $R$ is regular.
\item $R$ is noetherian and each acyclic complex of injective modules is contractible.
\item Countable direct sums of semi-injective complexes are semi-injective.
\item Arbitrary colimits of semi-injective complexes are semi-injective.
\end{enumerate}
\end{theorem}

\begin{proof}
Injective modules are semi-injective as complexes so conditions (3) and (4) imply that countable direct sums of injective modules are injective, and so the ring $R$ is noetherian, by \ref{bass}. Thus, in the remainder of the proof we assume $R$ is noetherian. The triangulated category $\kinj$ is then compactly generated, by \ref{kinj:cg}, and the localization functor $Q\col \kinj\to \dcat$ is compatible with coproducts; this fact will be used without further remark.

(1)$\iff$(2)  It suffices to verify that $R$ is regular if and only if the functor $Q$ is an equivalence, by Lemma~\ref{ikac=0}. In view of \ref{cg:functor}, the desired result is a consequence of \ref{kinj:cg} and \ref{dcat:regular}.

(2)$\implies$(4) Since $R$ is noetherian, a colimit of complexes of injectives is also a complex of injectives. Thus Proposition~\ref{id=hid} provides the desired conclusion.

(4)$\implies$(3) is clear.

(3)$\implies$(2) It is enough to prove that each complex $I$ of injective modules is semi-injective; see Proposition~\ref{id=hid}.  For each integer $n$ set $I(n)=I_{\les n}$; this is a subcomplex of $I$. Evidently, $I=\cup_{n\in\BN}I(n)$ so, with $\iota_{n}\col I(n)\subseteq I(n+1)$ the inclusion, there is an exact sequence of complexes of $R$-modules
\[
0\lra \bigoplus_{n\in\BN}I(n) \xra{\ \theta\ } \bigoplus_{n\in\BN}I(n)\lra I\lra 0\,,
\]
where $\theta(x_{n}) = (x_{n} - \iota_{n-1}(x_{n-1}))$. Each $I(n)$ is semi-injective, since it is a complex of injectives with $I(n)_{i}=0$ for $i>n$, so the hypothesis yields that the direct sums above are semi-injective. It follows that the complex $I$ is semi-injective as well.
\end{proof}

To put the next result in context, we recall that the semi-injective property does not localize; see \cite[6.5]{Neeman:1996a} and \cite{Chen/Iyengar:2009a} for counter-examples.

\begin{corollary}
\label{regular}
Let $R$ be a commutative regular ring and $U$ a multiplicatively closed subset of $R$. If $I$ is a complex of injective $R$-modules, then the complex $U^{-1}I$ is semi-injective over the ring $R$ and also over the ring $U^{-1}R$.
\end{corollary}

\begin{proof}
The ring $U^{-1}R$ is regular, by a result of Auslander, Buchsbaum, and Serre; see \cite[19.3]{Matsumura:1986a}. If $E$ is an injective $R$-module, then $U^{-1}E$ is injective over $R$ and over $U^{-1}R$. The  assertion is thus a direct consequence of Theorem~\ref{thm:inj}.
\end{proof}

Next we present a variation of Theorem~\ref{thm:inj} involving covers for complexes.

Following Enochs, Jenda, and Xu \cite{Enochs/Jenda/Xu:1996a}, we say that a graded-injective complex $I$ is \emph{minimal} if for each integer $n$ the inclusion $\ker(\dd_{n})\subseteq I_{n}$ is an essential extension.

Let $\eps\col E\to X$ be a morphism of complexes with $E$ acyclic. We say that $\eps$ is an \emph{acyclic cover} (or, as in \cite{Enochs/Jenda/Xu:1996a}, an \emph{exact cover}) of $X$ if each morphism $F\to X$ of complexes with $F$ acyclic factors uniquely through $\eps$. In~\cite[3.18]{Enochs/Jenda/Xu:1996a} it is proved that $\eps$ is an acyclic cover of $X$ if and only if $\eps$ is surjective and $\ker(\eps)$ is a minimal semi-injective complex. This result is used without comment in the proof below.

\begin{proposition}
\label{covers}
A ring $R$ is regular if and only if each direct sum of acyclic covers of complexes of $R$-modules is an acyclic cover.
\end{proposition}

\begin{proof}
Suppose $R$ is regular. Let $\eps_{\lambda}\col E_{\lambda}\to X_{\lambda}$, with $\lambda$ some index set, be a family of acyclic covers. The complexes $\ker(\eps_{\lambda})$ are semi-injective and minimal; hence
the complex $K=\oplus_{\lambda}\ker(\eps_{\lambda})$ is also semi-injective, by Theorem~\ref{thm:inj}, and
minimal, since minimality is preserved under direct sums. Setting $E=\oplus_{\lambda}E_{\lambda}$ and $X=\oplus_{\lambda}X_{\lambda}$, one thus obtains  an exact sequence of complexes
\[
0\to K \to E \to X  \to 0
\]
with $E$ acyclic and $K$ minimal semi-injective. Hence $E$ is an acyclic cover of $X$.

Assume now that each direct sum of acyclic covers is an acyclic cover. Given Theorem~\ref{thm:inj},  it suffices to verify if that $\{I_{\lambda}\}$ is an arbitrary family of semi-injective complexes, then $\oplus_{\lambda}I_{\lambda}$ is also semi-injective.

We may assume that each $I_{\lambda}$ is minimal. Indeed, any complex $I$ of injective $R$-modules is isomorphic to $I'\oplus I''$ with $I'$ minimal and $I''$ contractible; see \cite{Enochs/Jenda/Xu:1996a}
or \cite{Avramov/Foxby/Halperin:2008a}. When $I$ is itself semi-injective, so is $I''$ and hence it is homotopic to zero. Thus $I$ is homotopically equivalent to $I'$, which is minimal and semi-injective.

With $C_{\lambda}$ the mapping cone of the identity map of $I_{\lambda}$ one then obtains that the canonical morphism $\eps_{\lambda}\col C_{\lambda}\to \shift I_{\lambda}$ is an acyclic cover, by \cite[3.21]{Enochs/Jenda/Xu:1996a}. Thus our hypothesis implies that the morphism
\[
\eps=\bigoplus_{\lambda}\eps_{\lambda}\col
\bigoplus_{\lambda}C_{\lambda}\lra \bigoplus_{\lambda} \shift I_{\lambda}
\]
is an acyclic cover as well. Therefore the complex $\ker(\eps)$, that is to say, $\oplus_{\lambda} I_{\lambda}$, is semi-injective, as desired.
\end{proof}

\begin{remark}
Proposition~\ref{id=hid} and Theorem~\ref{thm:inj} can be generalized, using similar arguments, to locally noetherian Grothendieck categories with compactly generated derived categories. In particular, they yield an analogue of Theorem~\ref{ithm:inj} for quasi-coherent sheaves over a noetherian scheme.
\end{remark}

\section{Complexes of flat modules and of projective modules}
\label{Projectives}
We  present analogues of results in Section~\ref{Injectives} for complexes of flat modules and of projective modules. Many arguments are similar, so details are provided only when there are noteworthy differences.

The notions of graded-projective complexes and semi-projective complexes are obvious analogues of that of graded-injective complexes and semi-injective complexes; see ~\cite{Avramov/Foxby/Halperin:2008a}, or \cite{Avramov/Foxby:1991a} where they are called \#-projective and DG projective complexes, respectively. The invariants of interest are the \emph{projective dimension}:
\[
\pd_{R}M = \inf
\left\{n\in\BZ  \left|
\begin{gathered}
\text{$M\simeq P$ with $P$ semi-projective}\\
\text{and $P_{j}=0$ for $j>n$.}
\end{gathered}
\right.
\right\}
\]
Allowing $P$ above to any graded-projective yields an invariant denoted $\hpd_{R}M$.

The statement and proof of the next result parallel Proposition~\ref{id=hid}. Contractible complexes of projectives are the \emph{categorically projective} complexes of \cite{Avramov/Foxby/Halperin:2008a}.

\begin{proposition}
\label{pd=hpd}
Let $R$ be a ring. The following conditions are equivalent.
\begin{enumerate}[\quad\rm(1)]
\item An equality $\hpd_{R}M=\pd_{R}M $ holds for each complex $M$ of $R$-modules.
\item An equality $\hpd_{R}M=\pd_{R}M $ holds for each $R$-module $M$.
\item Each complex of projective $R$-modules is semi-projective.
\item Each acyclic complex of projective $R$-modules is contractible.  \qed
\end{enumerate}
\end{proposition}

Considering, in the same vein as before,  graded-flat complexes and semi-flat complexes generates invariants that we denote $\hfd_{R}M$ and  $\fd_{R}M$, respectively. A complex $F$ of $R$-modules is said to be \emph{categorically flat} if each module $F_{i}$ is flat and $\hh{M\otimes_{R}F}=0$ for each right $R$-module $M$. Once again the notions are from \cite{Avramov/Foxby:1991a} but terminology is from \cite{Avramov/Foxby/Halperin:2008a}. Categorically flat complexes have also been called `flat complexes'; see, for instance,  \cite{Enochs/Rozas:1998a}.

\begin{chunk}
\label{ch:flat}
Let $F$ be a complex of flat $R$-modules. The following conditions are equivalent:
\begin{enumerate}[\quad\rm(1)]
\item $F$ is categorically flat.
\item $\Hom RPF$ is acyclic for each complex $P$ of projective $R$-modules.
\item $F$ is acyclic and the $R$-module $\ker(\dd^{F}_{i})$ is flat for each $i$.
\end{enumerate}

Indeed, (1) and (3) are readily seen to be equivalent; the equivalence of (2) and (3) is due to Neeman~\cite[8.6]{Neeman:2008a}.
\end{chunk}

The next result can be proved along the same lines as Proposition~\ref{id=hid}.

\begin{proposition}
\label{fd=hfd}
Let $R$ be a ring. The following conditions are equivalent.
\begin{enumerate}[\quad\rm(1)]
\item An equality $\hfd_{R}M=\fd_{R}M $ holds for each complex $M$ of $R$-modules.
\item An equality $\hfd_{R}M=\fd_{R}M $ holds for each $R$-module $M$.
\item Each complex of flat $R$-modules is semi-flat.
\item Each acyclic complex of flat $R$-modules is categorically flat.  \qed
\end{enumerate}
\end{proposition}

We write $\kproj$ for the homotopy category of complexes of projective modules, viewed as a triangulated category,  and $\kpac$ for its full triangulated subcategory consisting of acyclic complexes; see ~\cite{Verdier:1996a}. The canonical localization functor $\kproj\to \dcat$ is again denoted $Q$.

The next result shows that conditions in Propositions ~\ref{pd=hpd} and \ref{fd=hfd} are equivalent.

\begin{proposition}
\label{kpac=0}
Let $R$ be a ring. The following conditions are equivalent.
\begin{enumerate}[\quad\rm(1)]
\item The localization functor $Q\col \kproj \to \dcat$ is an equivalence.
\item Each acyclic complex of projective $R$-modules is contractible.
\item Each acyclic complex of flat $R$-modules is categorically flat.
\end{enumerate}
\end{proposition}

\begin{proof}
(1)$\iff$(2) holds because the kernel of the functor $Q$ is $\kpac$.

(2)$\implies$(3) Let $F$ be an acyclic complex of flat $R$-modules. For each complex $P$ of projective modules
the complex $\Hom RPF$ is then acyclic, since $P$ is semi-projective by Proposition~\ref{pd=hpd}. Hence $F$ is categorically flat, by~\ref{ch:flat}.

(3)$\implies$(2) Let $P$ be an acyclic complex of projective $R$-modules. The complex $P$ is categorically flat, so $\Hom RPP$ is acyclic, by \ref{ch:flat}. Therefore $\idmap^{P}$ is homologous to zero; equivalently, $P$ is contractible.
\end{proof}

In what follows we write $\dcatfp$ for full subcategory of $\dcat$ consisting of complexes isomorphic to bounded complexes of finitely presented right $R$-modules. When $R$ is noetherian this coincides with the subcategory $\dcatf$.

The result below is due to J{\o}rgensen~\cite[3.2]{Jorgensen:2005a} under additional hypotheses on $R$;
the general case is contained in the work of Neeman~\cite[7.12, 7.14]{Neeman:2008a}.

\begin{chunk}
\label{kproj:cg}
When $\op R$ is coherent $\kproj$ is compactly generated, and the composition of functors $\Hom R-R\col \kproj\to \ck(\op R)$ and localization $\ck(\op R)\to \dcat[\op R]$ induces  an equivalence of categories:
\[
\kprojc \xra{\sim}\dcatfp[\op R]\,.
\]
\end{chunk}

We require only the following consequence of the result above:

\begin{lemma}
\label{lem:neeman}
Assume that the ring $\op R$ is coherent. The localization functor $Q\col \kproj \to \dcat$ induces an equivalence of categories $\kprojc\xra{\sim}\cd^{\ssc}(R)$ if and only if each complex in $\dcatfp[\op R]$ is perfect.
\end{lemma}

\begin{proof}
The equivalence in \ref{kproj:cg} implies that each complex in $\dcatfp[\op R]$ is perfect if and only if
each complex in $\kprojc$ is isomorphic to a bounded complexes of finitely generated projective modules;
that is to say,  when $Q$ induces an equivalence of categories $\kprojc\xra{\sim}\cd^{\ssc}(R)$, by \ref{dcat:cg}.
\end{proof}

The following characterization of coherent rings is due to Chase~\cite[2.1]{Chase:1960a}.

\begin{chunk}
\label{chase}
A ring $R$ is right coherent if and only if a product of flat $R$-modules is flat.
\end{chunk}

The theorem below and Propositions~\ref{pd=hpd} and \ref{fd=hfd} imply Theorem~\ref{ithm:prj}.

\begin{theorem}
\label{thm:prj}
Let $R$ be a ring. The following conditions are equivalent.
\begin{enumerate}[\quad\rm(1)]
\item $\op R$ is coherent and each complex in $\dcatfp[\op R]$ is perfect.
\item $\op R$ is coherent and each acyclic complex of projective modules is contractible.
\item Products of semi-flat complexes are semi-flat.
\end{enumerate}
\end{theorem}

\begin{proof}
Condition (3) implies that $\op R$ is coherent; this is by \ref{chase}. Thus, in the remainder of the proof we assume $\op R$ is coherent.

(1)$\iff$(2) The triangulated category $\kproj$ is compactly generated, by \ref{kproj:cg}. Given \ref{cg:functor}, the desired result now follows from  Proposition~\ref{kpac=0} and Lemma~\ref{lem:neeman}.

(2)$\implies$(3) By Proposition~\ref{kpac=0}, acyclic complexes of flat modules are categorically flat. It follows from \ref{chase} that a product of complexes of flat modules is a complex of flat modules, so the desired implication is a consequence of Proposition~\ref{fd=hfd}.

(3)$\implies$(2)
By Propositions~\ref{fd=hfd} and \ref{kpac=0}, it suffices to prove that if $F$ is a complex of flat $R$-modules, then it is semi-flat. Let $F(n)= F/F_{\sls - n}$ for each integer $n\geq 0$ and $\eps(n)\col F(n+1)\to F(n)$ the obvious surjection. The complex $F$ is the limit of the surjective system $\cdots \to F(n+1)\to F(n)\to \cdots$ so there is an exact sequence
\[
0\lra F\lra \prod_{n\in\BN}F(n)\xra{\ \nu\ } \prod_{n\in\BN}F(n)\lra 0\,,
\]
of complexes of flat $R$-modules, where $\nu(x_{n}) = (x_{n}-\eps(n)(x_{n+1}))$. Since each $F(n)$ is a complex of flat modules with $F(n)_{i}=0$ for $i<n$, it is semi-flat, and hence the complex $\prod_{n\in\BN}F(n)$ is also semi-flat, by hypothesis. The exact sequence above implies that $F$ is semi-flat, as desired.
\end{proof}

\begin{ack}
It is a pleasure to thank Luchezar Avramov for sharing his thoughts about the material discussed here, and the referee for detailed comments on previous versions of this article.
\end{ack}

\end{document}